\documentclass[a4paper,11pt,oneside]{book}
\usepackage[english]{certus}
\usepackage{parskip}
\usepackage{fancyhdr}
\title{CONSTRUCTION OF DIFFERENCE SCHEMES FOR NONLINEAR SINGULAR PERTURBED EQUATIONS BY
APPROXIMATION OF COEFFICIENTS}
\author{L.~V.~Rozanova\\[1em]{\em Scientific Advisor}\\A.~I.~Zadorin}
\date{Omsk -- 2000}
\begin{document}
\maketitle
\tableofcontents

\phantomsection
\addcontentsline{toc}{chapter}{Relevance}
\chapter*{Relevance\markboth{Relevance}{}}
Mathematical modeling of many physical processes such as
diffusion, viscosity of fluids and combustion involves differential equations
with small coefficients of higher derivatives. These may be small diffusion coefficients
for modeling the spreading of impurities, small coefficients of viscosity in fluid flow
simulation etc.

The difficulty with solving such problem is that if you set the small parameter at higher
derivatives to zero, the solution of the degenerate problem doesn't correctly approximate
the original problem, even if the small parameter approaches zero; the solution
of the original problem exhibits the emergency of a boundary layer. As a result, the
application of classical difference schemes for solving such equations produces
great inaccuracies. Therefore, numerical solution of differential equations with
small coefficients at higher derivatives demands special difference schemes exhibiting
uniform convergence with respect to the small parameters involved.

\paragraph{Current status of this problem\\}
The development of numerical methods for problems with exponential and power-law boundary
layers introduced among others the following methods: 
\begin{enumerate}
 \item densification of grids in boundary layers;
 \item scheme adjustment on the boundary layer component of solutions.
\end{enumerate}

The first approach includes the works of N.~S.~Bakhvalov, V.~D.~Liseikin, G.~I.~Shishkin,
V.~B.~Andreyev, R.~Vulanovich and other authors. In P.~S.~Bakhvalov's approach, the grid
nodes distribution is constructed in such a way that approximation inaccuracy on the
grid nodes of the boundary layers stays equal while outside of the boundary layer the mesh
remains uniform.

It is shown that application of such a grid yields second-order accuracy with respect to
the number of grid points. V.~D.~Liseikin purposes to perform change of variables
in such a way that derivatives up to some order are uniformly bounded. In terms of 
original variables, this is equivalent to mesh densification. G.~I.~Shishkin defined
an approach for building grids that are uniform both within the boundary layer and outside
of the boundary layer. In the works of G.~I.~Shishkin, V.~B.~Andreev, E.~A.~Savin,
N.~V.~Kopteva it is shown that on such a grid, a number of difference schemes (including
non-monotone scheme of central differences) are uniformly convergent. The significance of
this approach for partial differential equations is strengthened by the fact that in the case
of a parabolic boundary layer, as shown by G.~I.~Shishkin, there is no uniformly
convergent adjustment schemes on a uniform grid.

The second approach is presented in the works of A.~M.~Ilyin, G.~I.~Shishkin,
K.~V.~Emelyanov, D.~Miller, R.~Kellogg and others. The main idea of this approach is
isolating the boundary layer component of solutions and constructing a difference scheme
which is exact on the boundary layer function. The advantage of this approach is that it
does not impose restrictions on the mesh gauge, while the drawback is that the
boundary layer function must be expressed explicitly and the scheme must be specifically
adjusted to the function, which is not always possible.\\

In this article author investigates two nonlinear boundary value problems on a
finite interval, resulting in exponential and power-law boundary layers.

In the first chapter author constructs a first-order accurate finite difference scheme
for an exponential boundary layer problem. Author then introduces a certain grid where on
each interval the coefficients are replaced by constants. The second chapter
provides a construction of a first-order accurate difference scheme for a problem with
power-law boundary layer using a similar technique with additional mesh densification on
the boundary layer for scheme convergence.

\chapter{Nonlinear differential equation with exponential boundary layer}
The purpose of this chapter is to construct a difference scheme for a nonlinear equation
with exponential boundary layer, uniformly convergent in the small parameter. Since the
presence of a small parameter $\eps$ at the second-order derivative causes an exponential
growth of solutions near the boundaries of the interval, application of classical
difference schemes for solution of equations with small coefficient at the highest
derivative produces great inaccuracies. This follows from the fact that the accuracy of
known schemes is estimated as a product of some power of the mesh size $h$ and maximum
absolute value of some derivative of the solution, and derivatives of the solutions grow
unbounded in boundary layer with decreasing $\eps$. Thus the classical accuracy estimation
of difference schemes for solving equations with small parameter is not
acceptable.

Let us consider the boundary value problem:
\begin{equation}
 \left\{\begin{aligned}
 & T_\eps u = -\eps u'' + a(x)\, u' + g(u) = 0\\
 & u(0) = A,\quad R_\eps u = \eps u(L) + f(u(L)) = 0
 \end{aligned}\right.\label{eq:ra1}
\end{equation}
assuming that:
\begin{equation}
 \left\{\begin{aligned}
 & a(x) \ge \alpha, \quad \frac{\partial g}{\partial u} \ge -\beta, \quad \frac{\partial
f}{\partial u} \ge 0, \\
 & \alpha > 0,\ \beta > 0,\ \eps > 0, \ \alpha^2 - 4 \eps \beta \ge \gamma > 0.
\end{aligned}\right.
\end{equation}

Let's assume the function $g(s)$ to be twice continuously differentiable for all $s \in
\mathbb{R}$, $\eps \in (0, 1]$, and $a(x)$ continuously differentiable.

$C$ and $C_i$ will everywhere designate positive constants, not dependent on $\eps$ and
the mesh gauge.

\section{Analysis of the original problem}
Define an auxiliary linear operator $L$:
$$ Lu = -\eps u'' + a(x)\, u' + c(x)\,u $$
with boundary conditions:
$$ u(0), \quad Du = \eps u'(L) + d(x)\,u.$$
Let us find out whether for the operator $L$ the Principle of Maximum holds.

\begin{Lemma}
Assume $\exists \phi(x) \ge 0: L\phi(x) > 0, \ D\phi > 0 $. Then for the operator $L$
the Principle of Maximum holds, i.e. for an arbitrary twice differentiable function
$\psi(x)$ the condition:
$$ L\psi(x) \ge 0,\ D\psi \ge 0, \ \psi(0) \ge 0$$
implies that: $\psi(x) \ge 0$. \label{lem:l1}
\end{Lemma}
\begin{proof}
Assume $\psi(x)$ to be some differentiable function. Express $\psi(x)$ as the product:
$\psi(x) = v(x) \phi(x)$.\\
Assume there is $x_0 : \psi(x_0) < 0$.

Then $v(x_0) < 0$.

$$ v(0) \ge 0,\ D\psi = (v\phi)'(L) + dv\phi(L) \ge 0 $$
$$ v'\phi + v\phi' + dv\phi \ge 0.$$
$$ v'\phi + [\phi' + d\phi]v \ge 0.$$

Consequently,
$$ v'\phi + D\phi v \ge 0.$$

\begin{align}
 L\psi &= L(v\phi) = -\eps(v\phi)'' + a(v\phi)' + cv\phi \notag\\
       &= -\eps(v'\phi + v\phi')' + a(v'\phi + v\phi') + cv\phi \notag\\
       &= -\eps(v''\phi + 2 v' \phi' + v\phi'') + a(v'\phi + v\phi') + cv\phi \notag\\
       &= -\eps(v''\phi + 2 v' \phi') + v(-\eps \phi'' + a\phi + c\phi) + av'\phi \notag\\
       &= -\eps v''\phi + 2 \eps v' \phi' + av'\phi + vL\phi \notag\\
       &= -\eps\phi v'' - v'(2\eps\phi' - a\phi) + vL\phi \ge 0 \label{eq:longue}
\end{align}

Consider two cases:
\begin{enumerate}
 \item Let $v(L) \ge 0$. Then there is $\eta$, such that $v(\eta) < 0$ holds on the
interval $[0,L]$. Then $\eta$ is a minimum point for the function $v(x)$. At the point of
minimum first derivative equals zero, therefore, given the conditions of Lemma, in
\eqref{eq:longue} $\eps \phi v'' \ge 0$, $vL\phi > 0$, and the term $v'(2\eps \phi' -
a\phi)$ vanishes. So $L\psi < 0$, but from condition of the Lemma $L\psi \ge 0$,
a contradiction follows.
 \item Let $v(L) < 0$. Then due to the inequality $$ v' \phi + D\phi v \ge 0$$
 we obtain
 $$ v'(L) > 0$$
 and, consequently, there is a minimum point $\tilde{\eta}$. Arguing similarly to the
first case, we arrive at a contradiction.
\end{enumerate}
\end{proof}

Suppose that 
$$ d(x) \ge 0,\ c(x) \ge -\beta,\ \alpha^2 - 4\beta\eps \ge \gamma > 0,\ \beta > 0. $$
We show that for $L$ the Principle of Maximum holds. Define the function $\phi(x)$ as:
$$ \phi(x) = e^{\frac{\delta}{\eps}(x - L)}, \text{ where } \delta =
\frac{2\beta\eps}{\alpha}.$$

$$ D \phi(L) = \delta + d > 0.$$
\begin{align*}
L \phi(x) &=   \left(-\frac{\delta^2}{\eps} + a\frac{\delta}{\eps} + c\right) \phi(x)\\
          &\ge \left(-\frac{\delta^2}{\eps} + a\frac{\delta}{\eps} - \beta\right) \\
          &=   \frac{\beta}{\alpha^2}(\alpha^2 - 4\beta\eps)\,\phi(x) > 0.
\end{align*}

Thus,
$$ \phi > 0,\ D\phi > 0,\ L\phi > 0,$$
i.e. the Principle of Maximum holds.

Now we obtain an estimate for stability for the problem \eqref{eq:ra1}.
\begin{Lemma}
 Let $p(x)$ and $q(x)$ be two arbitrary functions. Then
 $$ \norm{p(x) - q(x)} \le c\norm{T_\eps p - T_\eps q} + c \abs{p(0) - q(0)} + c
\abs{R_\eps p - R_\eps q}.$$\label{lem:ll2}
\end{Lemma}
\begin{proof}
 Define $z = p - q$.\\
 Then with respect to $z$ we obtain the boundary problem:
 $$ Lz = -\eps z'' + a(x)\, z' + c(x)\, z = T_\eps p - T_\eps q,\text{ where} $$
 $$ c(x) = \frac{g(p) - g(q)}{p - q} \ge -\beta. $$
 $$ z(0) = p(0) - q(0).$$

 \begin{align*}
 Dz &= \eps z'(L) + dz(L) = \eps z' (L) + \frac{f(p(L)) - f(q(L))}{p(L) - q(L)} z(L) \\
    &= \eps z'(L) + f'(\theta) = R_\eps p - R_\eps q.
 \end{align*}
 Introduce the function
 $$ \psi(x) = c_1 e^{\frac{2\beta x}{\alpha}} \norm{T_\eps p - T_\eps q} + c_2
e^{\frac{\alpha}{\eps}(x - l)} \norm{R_\eps p - R_\eps q} + e^{\frac{2\beta x}{\alpha}}
\abs{p(0) - q(0)} \pm z(x).$$
 It is easy to verify that:
 \begin{align*}
  L e^{\frac{\alpha(x - L)}{2\eps}} & \ge \left(-\frac{\alpha^2}{4\eps} +
\frac{\alpha^2}{2\eps} + c(x)\right)\, e^{\alpha(x - L)}{2\eps} \ge
\left(-\frac{\alpha^2}{4\eps} - \beta\right)\, e^{\alpha(x - L)}{2\eps}\\
 &= \left(-\frac{\alpha^2 - 4\beta \eps}{4\eps}\right)\, e^{\alpha(x - L)}{2\eps} \ge
\gamma e^{\alpha(x - L)}{2\eps} > 0.
 \end{align*}

\begin{align*}
  L e^{\frac{2\beta x}{\alpha}} & \ge \left(-\eps \frac{4\beta^2}{\alpha^2} +
2\beta + c(x)\right)\, e^{2\beta x}{\alpha} \ge
\left(-\eps \frac{4 \beta^2}{\alpha^2} + \beta\right)\, e^{2\beta x}{\alpha}\\
 &= \beta e^{2\beta x}{\alpha} \left(\frac{-4 \eps \beta +
\alpha^2}{\alpha^2}\right) \ge \frac{\beta \gamma}{\alpha^2} e^{2\beta x}{\alpha}
> 0.
 \end{align*}
 
Consequently,
$$ L\psi \ge c_1 \frac{\beta \gamma}{\alpha^2} e^{2\beta x}{\alpha} \norm{T_\eps p -
T_\eps q} - \norm{T_\eps p - T_\eps q} \ge 0, $$
where $c_1 = \frac{\alpha^2}{\beta \gamma}$.

$$ D e^{\frac{2\beta x}{\alpha}} = \left( \eps \frac{2\beta}{\alpha} + d\right)
e^{\frac{2\beta x}{\alpha}} \ge 0,$$
$$ D e^{\frac{\alpha(x - L)}{2\eps}} = \alpha + d > 0.$$

Consequently,
$$ L\psi \ge (\alpha + d) c_2 \abs{R_\eps p - R_\eps q} - \abs{R_\eps p - R_\eps q} \ge 0,
$$
where $c_2 = \frac{1}{\alpha + d}$.

Thus, defining
$$\psi(x) = \frac{\beta \gamma}{\alpha^2} e^{\frac{2\beta x}{\alpha}} \norm{T_\eps p -
T_\eps q} + \frac{1}{a + d} e^{\frac{\alpha{x - L}}{\eps}} \pm z(x),$$
we obtain $\psi(0) \ge 0,\ D\psi \ge 0,\ L\psi \ge 0$,
hence by Lemma~\ref{lem:l1} $\psi(x) \ge 0$.
\end{proof}

\begin{Corr}
Lemma~\ref{lem:ll2} implies uniqueness and boundedness of problem~\eqref{eq:ra1}.
\end{Corr}
Proof of the boundedness:\\
Let $p(x) = u(x),\ q(x) = 0$, then by Lemma~\ref{lem:ll2},
$$ \abs{u(x)} \le \frac{\alpha^2}{\beta \gamma} e^{\frac{2\beta x}{\alpha}} \abs{g(0)} +
\frac{1}{a + d} e^{\frac{\alpha(x - L)}{\eps}} \abs{f(0)} + e^{\frac{2\beta x}{\alpha}}
\abs{u(0)}.$$

The uniqueness of this solution is obvious.

\section{Estimate of the derivative}
We obtain an estimate of the derivative problem~\eqref{eq:ra1}.
\begin{Lemma}
 $$ \abs{u'(x)} \le c\left[1 + \frac{1}{\eps} e^{\frac{a(x - L)}{\eps}}\right].$$
\end{Lemma}
\begin{proof}
 We express the equation~\eqref{eq:ra1} in the form:
 $$ \left(-\eps u' \exp\left[\int_x^L \frac{a(t)}{\eps}\,dt \right] \right)' + g(u(x))
\exp\left[\int_x^L \frac{a(t)}{\eps}\,dt \right] = 0$$
Integrating from $x$ to $L$:
$$ -\eps u'(L) + \eps u'(x) \exp\left[\int_x^L \frac{a(t)}{\eps}\,dt \right] + \int_x^L
g(u(s)) \exp\left[\int_x^L \frac{a(t)}{\eps}\,dt \right]\, ds = 0$$
$$ u'(x) = u'(L)\exp\left[-\int_x^L \frac{a(t)}{\eps}\,dt \right] - \frac{1}{\eps}
\int_x^L g(u(s)) \exp\left[-\int_x^L \frac{a(t)}{\eps}\,dt \right]\, ds $$

It is clear that:
\begin{align*}
 \frac{1}{\eps} \int_x^L g(u(s)) \exp\left[-\int_x^L \frac{a(t)}{\eps}\,dt \right]\, ds
 &\le \frac{c}{\eps} \int_x^L \exp\left[-\int_x^L \frac{a(t)}{\eps}\,dt \right]\, ds\\
 &=   \frac{c}{\eps} \int_x^L e^{\frac{\alpha}{\eps}(x - L)}\, ds\\
 &=   \frac{c}{\eps} \left[1 - e^{\frac{\alpha}{\eps}(x - L)}\right]\\
 &\le \frac{c}{\alpha} = c_1\\
\end{align*}

From the boundary condition we conclude that
$$ \abs{u'(x)} \le \abs{\frac{f(u(L))}{\eps}} \le \frac{c}{\eps}, $$
consequently,
$$ \abs{u'(x)} \le c + \frac{c}{\eps} e^{\frac{\alpha}{\eps}(x - L)}. $$
\end{proof}

The obtained estimate of the derivative characterizes the boundary layer.

\section{Construction of a difference scheme}
Let us introduce a non-uniform grid and replace the coefficients by constants in each
interval of that grid. That allows us to write the solution explicitly.
Matching derivatives on adjacent intervals will lead to a difference scheme.\\

So, let $\Delta_n = [x_{n - 1}, x_n]$.\\

Now we turn to a problem with piecewise constant coefficients:
\begin{equation*}
 \left\{\begin{aligned}
 & \eps V'' - \tilde{a}V' + \tilde{g}(V) = 0,\\
 & V(0) = A,\quad R_\eps V = 0,
 \end{aligned}\right.
\end{equation*}
where $\tilde{a} = a_n = a(x_{n - 1})$, $\tilde{g}(V(x)) = g_n = g(V(x_{n - 1}))$, at
$x \in \Delta_n$.\\

On an arbitrary interval $\Delta_n$ we have:
\begin{equation*}
 \left\{\begin{aligned}
 & \eps V'' - a_n V' + g_n = 0,\quad x \in \Delta_n\\
 & V(x_{n - 1}) = V_{n - 1}^h,\quad V(x_{n}) = V_n^h,
 \end{aligned}\right.
\end{equation*}

where $\{V_n^h\}$ is not yet defined.\\
The solution on the interval $\Delta_n$ has the form:
\begin{equation}
 V(x) = c_1 + c_2 e^{\frac{a_n}{\eps} (x - x_n)} + \frac{g_n}{a_n}x.
\end{equation}

Let $h_n$ be the grid size: $h_n = x_n - x_{n-1}$.\\
Define $c_2$ from the boundary conditions:
\begin{equation*}
 \left\{\begin{aligned}
 & c_1 + c_2 e^{-\frac{a_n}{\eps}h_n} + \frac{g_n}{a_n} x_{n-1} = V_{n-1}^h,\\
 & c_1 + c_2 + \frac{g_n}{a_n} x_n = V_n^h.
 \end{aligned}\right.
\end{equation*}

Let us find $c_2$:
\begin{gather*}
 V_{n-1}^h - c_2 e^{-\frac{a_n h_n}{\eps}} - \frac{g_n}{a_n} x_{n-1} = V_n^h - c_2 -
 \frac{g_n}{a_n} x_n, \\
 V_{n-1}^h - V_n^h + \frac{g_n}{a_n} h_n = c_2 \left( e^{-\frac{a_n h_n}{\eps}} -
1\right),\\
 c_2 = \left.\left(V_{n-1}^h - V_n^h + \frac{g_n}{a_n} h_n\right) \middle/ \left(
e^{-\frac{a_n h_n}{\eps}} - 1\right)\right..
\end{gather*}

Then let us find $V'(x)$:
$$ V'(x) = \frac{a_n}{\eps} c_2 e^{\frac{a_n}{\eps}(x - x_n)} + \frac{g_n}{a_n}. $$

To ensure that the solution is continuously differentiable, the derivatives of
solutions must be matched on the adjacent the intervals. This requires
$$ \lim_{x \rightarrow x_{n - 0}} V'(x) = \lim_{x \rightarrow x_{n + 0}} V'(x) $$
$$ \frac{g_n}{a_n} + \frac{a_n}{\eps} c_2 = \frac{g_{n+1}}{a_{n+1}} + \frac{a_{n+1}}{\eps}
\tilde{c}_2 e^{-\frac{a_{n+1} h_{n+1}}{\eps}}, $$
$$ \frac{g_n}{a_n} + \frac{a_n}{\eps}\left.\left(V_{n-1}^h - V_n^h + \frac{g_n}{a_n}
h_n\right) \middle/ \left( e^{-\frac{a_n h_n}{\eps}} - 1\right)\right. = $$
$$ = \frac{g_{n+1}}{a_{n+1}} + \frac{a_{n+1}}{\eps} \left.\left(V_{n}^h - V_{n+1}^h +
\frac{g_{n+1}}{a_{n+1}} h_{n+1}\right) \middle/ \left( e^{-\frac{a_{n+1} h_{n+1}}{\eps}} -
1\right)\right.$$

$$ \frac{g_n}{a_n} + \left.\left(\frac{V_{n}^h - V_{n-1}^h}{h_n} -
\frac{g_n}{a_n}\right) \middle/ \frac{\eps}{h_n a_n} \left(1 - e^{-\frac{a_n h_n}{\eps}}
-\right)\right. = $$
$$ = \frac{g_{n+1}}{a_{n+1}} + \left.\left(\frac{V_{n+1}^h - V_n^h}{h_{n+1}} -
\frac{g_{n+1}}{a_{n+1}}\right) e^{-\frac{a_{n+1} h_{n+1}}{\eps}} \middle/
\frac{\eps}{h_{n+1} a_{n+1}} \left(1 - e^{-\frac{a_{n+1} h_{n+1}}{\eps}}
-\right)\right..$$

Let $s_n = \frac{\eps}{h_n a_n} \left(1 - e^{-\frac{a_{n}
h_{n}}{\eps}}\right)$.\\

Then the scheme takes the form:
$$ \frac{g_n}{a_n} + \left.\left(\frac{V_{n}^h - V_{n-1}^h}{h_n} -
\frac{g_n}{a_n}\right) \middle/s_n \right. = $$
$$ = \frac{g_{n+1}}{a_{n+1}} + \left.\left(\frac{V_{n+1}^h - V_n^h}{h_{n+1}} -
\frac{g_{n+1}}{a_{n+1}}\right) e^{-\frac{a_{n+1} h_{n+1}}{\eps}} \middle/s_{n+1}\right..$$

Rewrite the scheme in the form:
$$ \frac{g_n}{a_n} + \frac{V_{n}^h - V_{n-1}^h}{s_n h_n} - \frac{g_n}{s_n a_n} =
\frac{g_{n+1}}{a_{n+1}} + \left(\frac{V_{n+1}^h - V_{n}^h}{s_{n+1} h_{n+1}} -
\frac{g_{n+1}}{s_{n+1} a_{n+1}} \right) e^{-\frac{a_{n+1} h_{n+1}}{\eps}}.$$So,
\begin{align*}
 \frac{V_{n+1}^h - V_{n}^h}{s_{n+1} h_{n+1}}e^{-\frac{a_{n+1} h_{n+1}}{\eps}} -
 \frac{V_{n}^h - V_{n-1}^h}{s_n h_n} &=
\begin{multlined}[t]
\frac{g_{n+1}}{s_{n+1} a_{n+1}} e^{-\frac{a_{n+1} h_{n+1}}{\eps}} - \\
- \frac{g_n}{s_n a_n} + \frac{g_n}{a_n} - \frac{g_{n+1}}{a_{n+1}}.
\end{multlined}
\end{align*}

Supplement this scheme by an approximation of the boundary conditions:
$$ V_0^h = A,\ \eps V' + f(V(L)) = 0.$$

In the last interval $V'(x) = c^(N)\frac{a_N}{\eps} e^{\frac{a_N}{\eps}(x-L)}+
\frac{g_N}{a_N}$,
therefore, boundary conditions take the following form:
$$ V_0^h = A,$$
$$ c_2^N a_N + \frac{g_N\eps}{a_N} + f(V_N^h) = 0, \text{ i.e.}$$
$$ \frac{(V_N^h - V_{N-1}^h) a_N}{1 - e^{-\frac{a_N h_N}{\eps}}} - \frac{g_N h_N}{1 -
e^{-\frac{a_N h_N}{\eps}}} + \frac{g_N \eps}{a_N} + f(V_N^h) = 0$$
So, as a result of the difference scheme has the form:
\begin{align}
 &\frac{V_{n+1}^h - V_{n}^h}{s_{n+1} h_{n+1}}e^{-\frac{a_{n+1} h_{n+1}}{\eps}} -
 \frac{V_{n}^h - V_{n-1}^h}{s_n h_n} =
\begin{multlined}[t]
\frac{g_{n+1}}{s_{n+1} a_{n+1}} e^{-\frac{a_{n+1} h_{n+1}}{\eps}} - \\
- \frac{g_n}{s_n a_n} + \frac{g_n}{a_n} - \frac{g_{n+1}}{a_{n+1}}.
\end{multlined}\notag\\
& V_0^h = A,\ n = 1,2,\dots,(N-1)\notag\\
&\frac{(V_N^h - V_{N-1}^h) a_N}{1 - e^{-\frac{a_N h_N}{\eps}}} - \frac{g_N h_N}{1 -
e^{-\frac{a_N h_N}{\eps}}} + \frac{g_N \eps}{a_N} + f(V_N^h) = 0\label{eq:e13}
\end{align}

Due to $g_n = g(V_{n-1}^h)$, the scheme (1.3) is nonlinear and can be linearized by the
method of iterations \cite{OR}, and solved by sweeping at each iteration.\cite{Sam}.

\section{Evaluation of convergence}
\begin{Thm}
 Let $u(x)$ be solution of the problem~\eqref{eq:ra1}, $u^h$ the solution of the
scheme~\eqref{eq:e13}.
 Then there exists a constant $C$, independent of $\eps$, such that
 $$ \max_n \abs{u_n^h - u(x_n)} \le C \cdot \max_n h_n.$$
\end{Thm}
\begin{proof}
 Rewrite the problem~\eqref{eq:ra1} and auxiliary problem:
 \begin{align}
  T_\eps u &= -\eps u'' + a(x)u' + g(u), & u(0) &= A & \eps u'(L) + f(u(L)) &= 0,\notag\\
  \tilde{T}_\eps v &= -\eps v'' + a(x)v' + \tilde{g}(v), & v(0) &= A & \eps v'(L) +
f(v(L)) &= 0.\label{eq:a14}\\
 \end{align}
 We use the fact that the scheme~\eqref{eq:e13} is exact for the problem~\eqref{eq:a14}.
 Therefore it is sufficient to estimate the proximity of these problems.

 Let $z = u - v$.

 Write the problem on $z$:
 \begin{equation*}
 \left\{\begin{aligned}
 & -\eps z'' + az' + g(u) - \tilde{g}(v) = 0,\\
 & z(0) = 0,\ \eps z'(L) + f(u(L)) - f(v(L)) = 0.
 \end{aligned}\right.
 \end{equation*}

 Rewrite the problem in the form:
 $$-\eps z'' + az' + g(u) - \tilde{g}(v) + \tilde{g}(u) - \tilde{g}(u) = 0 $$
 $$-\eps z'' + az' + \frac{\tilde{g}(u) - \tilde{g}(v)}{u - v}\cdot z = \tilde{g}(u) -
g(u) $$

 For $z(x)$ we obtain the boundary problem:
 \begin{equation*}
 \left\{\begin{aligned}
 & Lz = -\eps z'' + a(x)z' + g_u'(\theta)z = \tilde{g}(u) - g(u),\\
 & z(0) = 0,\ Dz = \eps z'(L) + f_u'(\theta)z(L) = 0.
 \end{aligned}\right.
 \end{equation*}

 Estimate $g(u) - \tilde{g}(u)$.\\
 For $x \in \Delta_n$ we have:
 $$ \tilde{a}(x) = a(x_{n-1}),\ \tilde{g}(u) = g(u(x_{n-1})).$$

 We obtain:
 \begin{align*}
 \abs{g(u(x)) - \tilde{g}(u(x))} &= \abs{g(u(x_n)) - g(u(x))}\\
  &\le \abs{\frac{\partial g}{\partial u}} \abs{u(x_n) - u(x)}\\
  &\le C\abs{u(x_n) - u(x)} \le C \abs{\int_{x_{n-1}}^x u'(s)\, ds} \\
  &\le C \int_{x_{n-1}}^x \abs{u'(s)}\, ds \\
  &\le \tilde{C} \int_{x_{n-1}}^x \left[1 + \frac{1}{\eps} e^{\frac{a(s -
L)}{\eps}}\right]\, ds \le Ch_n + \frac{Ch_n}{\eps} e^{\frac{a(x-L)}{2\eps}}
 \end{align*}

 Let $h = \max_n h_n$.\\
 Define function $\psi(x)$:
 $$ \psi(x) = \tilde{C}h \left[e^{\frac{a(x-L)}{2\eps}} + e^{\frac{2\beta
x}{\alpha}}\right] \pm z(x)$$

 Then:
 $$ \psi(0) \ge 0,\ D\psi \ge 0, L\psi \ge \tilde{C}h \frac{C_1}{\eps}
 e^{\frac{a(x-L)}{2\eps}} + \tilde{C}\cdot C_2 h - Ch_n - \frac{Ch_n}{\eps}
 e^{\frac{a(x - L)}{2\eps}} \ge 0$$
 if $\tilde{C} \le \frac{C}{C_2}$, $\tilde{C} \le \frac{C}{C_1}$.

 By the Principle of Maximum $\psi(x) \ge 0$.\\
 Consequently, $\abs{z(x)} \le Ch$.
\end{proof}

Now we have constructed a difference scheme for a nonlinear second-order
equation with exponential boundary layer and proved its uniform convergence in the small
parameter.

\chapter{Nonlinear differential equations with power-law boundary layer}
Consider the following boundary problem:
\begin{equation}
 \left\{\begin{aligned}
 & (\eps + x)^2\, u'' - f(x,u) = 0,\\
 & u(0) = A,\quad u(1) = B,
 \end{aligned}\right.
 \label{eq:ra2}
\end{equation}
under asumption that $\frac{\partial f}{\partial u} \ge \alpha > 0$, $\eps \in (0,1]$
and the function $f$ is continuously differentiable in its arguments.

Let's investigate the problem \eqref{eq:ra2}. We will prove that the solution contains
exponential boundary layer, construct a difference scheme and prove its convergence.

$C$ and $C_i$ will everywhere designate the positive constants independent of
$\eps$ and the mesh gauge.

\section{Boundedness of solutions}
We prove the boundedness for the solution of the problem \eqref{eq:ra2}.
\begin{Lemma}
 $$\norm{u(x)} \le \frac{1}{\alpha} \, \norm{f(x, 0)}, $$
 where $\norm{u(x)} = \max_a \abs{u(x)}$.
\end{Lemma}
\begin{proof}
 Write the equation \eqref{eq:ra2} in the following form:
 \begin{align*}
  Lu &= (\eps + x)^2\, u'' - \frac{f(x, u) - f(u, 0)}{u}\, u = f(x, 0)\\
  Lu &= (\eps + x)^2\, u'' - c(x)\, u = f(x, 0),
 \end{align*}
 where $c(x) \ge \alpha > 0$.

 Introduce the function $\psi(x) = \frac{\norm{f(x, 0)}}{\alpha} \pm u(x)$.
 $$ L\psi(x) \le 0, \psi(0) \ge 0, \psi(1) \ge 0, $$
 Hence, according to the Principle of Maximum $\psi(x) \ge 0$.\\
 $\frac{\norm{f(x, 0)}}{\alpha} \pm u(x) \ge 0$, consequently,
 $\abs{u(x)} \le \frac{1}{\alpha} \, \norm{f(x, 0)}$.
\end{proof}

\section{Estimation of the derivative}
\begin{Lemma}
 $$\abs{u'(x)} \le \frac{C}{\eps + x}.$$ \label{lem:l2}
\end{Lemma}
\begin{proof}
 From the equation \eqref{eq:ra2} we conclude that
 \begin{align}
 \abs{u''(x)} & \le \frac{C}{(\eps + x)^2} \notag \\
 u''(x) & = \frac{f(x,u)}{(\eps + x)^2}. \label{eq:e2}
 \end{align}
 Integrate \eqref{eq:e2} from $\xi$ to $x$:
 $$ u'(x) - u'(\xi) = \int_\xi^x \frac{f(x, u(x))}{(\eps + x)^2}\, dx, $$
 consequently,
 $$ \abs{u'(x) - u'(\xi)} = \abs{\int_\xi^x \frac{f(x, u(x))}{(\eps + x)^2}\, dx}. $$

 \begin{enumerate}
  \item Let $\xi \le x$. Then
  $$ \abs{u'(x) - u'(\xi)} \le \left[\frac{-C}{\eps + x}\right]_\xi^x = C
\left[\frac{1}{\eps + \xi} - \frac{1}{\eps + x} \right]$$
  \item Let $\xi > x$. Then
  $$ \abs{u'(x) - u'(\xi)} = \int_\xi^x \frac{C}{(\eps + x)^2}\, dx, $$
  consequently,
  $$ \abs{u'(x) - u'(\xi)} \le -C \left[\frac{1}{\eps + \xi} - \frac{1}{\eps + x}
\right]$$
 \end{enumerate}
 Thus for any $\xi$
 $$ \abs{u'(x) - u'(\xi)} \le C \abs{\frac{1}{\eps + \xi} - \frac{1}{\eps + x}}.$$
 Let $\eps$ be so that:
 $$ \text{} \frac{1}{2} u'(\xi) = u(1) -
u\left(\frac{1}{2}\right). $$
 Then
 $$ \abs{u'(\xi)} \le C_0 \text{ and } \frac{1}{2} \le \xi \le 1.$$
 Since
 $$ \abs{u'(\xi)} \le C_0 + C \abs{\frac{1}{\eps + \xi} - \frac{1}{\eps + x}},$$
 then
 $$ \abs{u'(\xi)} \le \frac{C_1}{\eps + x},$$
 We took into account that
 $$\frac{1}{\eps + \xi} \le \frac{1}{\eps + \frac{1}{2}} \le \frac{1}{\frac{1}{2}} = 2$$
 $$ 0 \le \frac{1}{\eps + \xi} \le 2,$$
 $$ \abs{u'(\xi)} \le C_0 + 2 + \frac{C}{\eps + x} \le \frac{C_1}{\eps + x},$$
 Hence the lemma follows.
\end{proof}

According to Lemma~\ref{lem:l2}, at the boundary $x = 0$ we have the power boundary layer.

\section{Construction of the scheme}
Let's proceed to the problem with piecewise constant coefficients.
\begin{equation}
 \left\{\begin{aligned}
 & (\eps + x)^2\, v'' - \tilde{f}(x,v) = 0,\\
 & v(0) = A,\quad v(1) = B,
 \end{aligned}\right.
 \label{eq:ra4}
\end{equation}
 where $\tilde{f}(x, u) = f_n = f(x_{n-1}, V(x_{n_1}))$, at $x \in \Delta_n = [x_{n-1},
x_n]$.\\[1ex]
To construct the scheme of equation~\eqref{eq:ra4} let's express $v(x)$ as:
\begin{align*}
 v''(x) &= \frac{f_n}{(\eps + x)^2},\\
 v'(x)  &= -\frac{f_n}{\eps + x} + C_1,\\
 v(x)   &= -f_n \cdot \ln(\eps + x) + C_1 x + C_2
\end{align*}

Let's label $v(x_{n - 1}) = V_{n-1}^h$, $v(x_n) = V_n^h$.

Let's find $C_1$ from the boundary conditions:
\begin{equation*}
 \left\{\begin{aligned}
 & - f_n \cdot \ln(\eps + x_{n-1}) + C_1 x_{n-1} + C_2 = V_{n-1}^h,\\
 & - f_n \cdot \ln(\eps + x_n) + C_1 x_n + C_2 = V_n^h,
 \end{aligned}\right. 
\end{equation*}
$$ C_1 h_n - f_n \, \ln \frac{\eps + x_n}{\eps + x_{n - 1}} = V_n^h - V_{n-1}^h $$
Then,
$$ C_1 = \frac{f_n}{h_n}\, \ln \frac{\eps + x_n}{\eps + x_{n - 1}} +
\frac{V_n^h - V_{n-1}^h}{h_n}$$

We match the derivatives of solutions to the ends of adjacent intervals:
\begin{align*}
 \lim_{x \rightarrow x_{n - 0}} V'(x) &= \lim_{x \rightarrow x_{n + 0}} V'(x)\\
 -\frac{f_n}{\eps + x_n} + C_1 &= -\frac{f_{n + 1}}{\eps + x_{n + 1}} + \tilde{C}_1\\
 \begin{multlined}[b]
  -\frac{f_n}{\eps + x_n} + \frac{f_n}{h_n}\times \\ \times \ln\left(\frac{\eps +
x_n}{\eps +
x_{n - 1}}\right) + \\ + \frac{V_n^h - V_{n-1}^h}{h_n}
 \end{multlined}
 &= \begin{multlined}[t]
     -\frac{f_{n + 1}}{\eps + x_{n + 1}} +
\frac{f_{n+1}}{h_{n+1}} \times \\
\times \ln\left(\frac{\eps + x_{n+1}}{\eps + x_n}\right) + \frac{V_{n+1}^h -
V_n^h}{h_{n+1}}
    \end{multlined}
\end{align*}

Thus we obtain:
$$ \frac{V_n^h - V_{n-1}^h}{h_n} - \frac{V_{n+1}^h - V_n^h}{h_{n+1}}= $$
\begin{equation}
 = \frac{f_n}{\eps + x_n} - \frac{f_n}{h_n}\, \ln \frac{\eps + x_n}{\eps + x_{n - 1}}
-\frac{f_{n + 1}}{\eps + x_{n + 1}} + \frac{f_{n+1}}{h_{n+1}}\, \ln \frac{\eps
+ x_{n+1}}{\eps + x_n} \label{eq:scheme}
\end{equation}
$$ V_0^h = A,\ V_N^h = B, \ n = 1, 2, 3,\dots, (N-1);\ f_n = f(x_{n-1}, V_{n-1}^h)$$

\section{Justification of convergence}
Define the mesh $\Omega$ which densifies in the boundary layer in such a way that
 $\ln\left(\frac{\eps + x_n}{\eps + x_{n-1}}\right) = \frac{C}{N}$.

Let $x_n = \lambda\left(\frac{n}{N}\right)$, then
$$ \ln\left(\eps + \lambda\left(\frac{n}{N}\right)\right) - \ln\left(\eps +
\lambda\left(\frac{n - 1}{N}\right)\right) = \frac{C}{N}.$$

Consequently,
$$ \frac{\lambda'\left(\frac{n}{N}\right)}{\eps + \lambda\left(\frac{n}{N}\right)} \cdot
\frac{1}{N} = \frac{C}{N},$$

That is 
$$ \frac{\lambda'\left(\frac{n}{N}\right)}{\eps + \lambda\left(\frac{n}{N}\right)} =
C,\ \lambda(0) = 0, \lambda(1) = 1,$$

Consequently,
$$ \lambda(t) = \eps\left[\left(1 + \frac{1}{\eps}\right)^t - 1\right].$$

Thus,
$$ \Omega = \left\{x_n = \eps\left[\left(1 +
\frac{1}{\eps}\right)^{\frac{n}{N}} - 1\right], \ n = 0,1,\dots,N\right\}.$$

\begin{Thm}
 $$ \norm{V^h - [u]_\Omega} \le \frac{C}{N}\, \ln\left(1 + \frac{1}{\eps}\right) $$
\end{Thm}
\begin{proof}
 Let $z = u - v$, where $u$ is the solution of problem \eqref{eq:ra2}. Since the scheme
\eqref{eq:scheme} is exact for the problem \eqref{eq:ra4} by its construction, it
is sufficient to estimate $\norm{v - v}$.

 Problem for $z$ has the form:
 \begin{equation}
 \left\{\begin{aligned}
 & (\eps + x)^2\, z'' - f'(\theta(x))z = f(x,v) - \tilde{f}(x,v),\\
 & z(0) = 0,\quad z(1) = 0,
 \end{aligned}\right.
 \label{eq:e6}
 \end{equation}

 We estimate the right side:
 \begin{multline}
  \abs{f(x,v) - \tilde{f}(x,v)} = \abs{f(x,v) - f(x_{n-1},v_{n-1})} \le \\
  \le C\abs{x - x_{n-1}} + C_1\abs{v(x) - v(x_{n-1})} \le Ch_n + C_1 \int_{x_{n-1}}^{x_n}
\abs{v'(s)}\, ds.
 \end{multline}

 Since $\abs{v'(x)} \le \frac{C_2}{x + \eps}$,
 $$ \abs{f(x,v) - \tilde{f}(x,v)} \le Ch_n + C_1 \int_{x_{n-1}}^{x_n}
\frac{C_2}{x + \eps}\, ds \le Ch_n + C_1 C_2 \ln\left(\frac{\eps + x_n}{\eps +
x_{n-1}}\right)$$

In accordance with the construction of the mesh:
$$ \abs{f(x,v) - \tilde{f}(x,v)} \le Ch_n + \frac{C_3}{N} \ln \left( 1 +
\frac{1}{\eps}\right) $$

We estimate $h_n$.
\begin{enumerate}
 \item From
 $$ x_n = \eps \left[ \left(1 + \frac{1}{\eps}\right)^{\frac{n}{N}} - 1\right] $$
 we conclude that
 $$ h_N = x_N - x_{N-1} = 1 - \eps \left[ \left(1 + \frac{1}{\eps}\right)^{\frac{n}{N}} -
1\right] \le C_5 \frac{\ln\left(1 + \frac{1}{\eps}\right)}{N} $$
 \item $h_n$ increases.
 Therefore it is sufficient to estimate $h_N$.\\
 Then $\abs{f(x,v) - \tilde{f}(x,v)} \le \frac{C}{N} \ln\left(1 +
\frac{1}{\eps}\right)$.\\
 Let $\psi(x) = \frac{C}{N} \ln\left(1 + \frac{1}{\eps}\right) \pm z(x)$.
 $$ \psi(0) \ge 0, \psi(1) \ge 0, L\psi(x) \le 0.$$
 Consequently,
 $$ \abs{z(x)} \le \frac{1}{\alpha} \cdot \frac{C}{N} \ln\left(1 + \frac{1}{\eps}\right).
$$
\end{enumerate}

 Thus,
 $$ \abs{u(x) - v(x)} \le \frac{C}{N} \ln\left(1 + \frac{1}{\eps}\right). $$

 Since $[V]_\Omega$ coincides with $V^h$, we obtain the theorem.
\end{proof}

\phantomsection
\addcontentsline{toc}{chapter}{Results}
\chapter*{Results\markboth{Results}{}}
At the result of the work, a first-order accurate difference scheme for a nonlinear
equation with exponential boundary layer and nonlinear boundary conditions is constructed
and its uniform convergence with respect to the small parameter is proven; a first-order
accurate difference scheme for a nonlinear equation with power-law boundary layer and
boundary conditions of first kind is constructed on a special mesh, and its convergence is
proven as well.

\phantomsection

\end{document}